\documentclass[reqno]{amsart}
\usepackage{graphics}
\numberwithin{equation}{section}

\newcommand{\iu}{\int_0^1}
\newcommand{\il}{\int_0^\infty}
\newcommand{\atn}{{\mathrm{tan}}^{-1}\,}
\newcommand{\D}[2]{\left(\frac{d}{d#1}\right)^{#2}}
\newcommand{\R}{\mathbf{R}}
\newcommand{\Cin}{\mathrm{Cin}}

\begin{document}

\title[Estimating Derivatives]{Using Integral Transforms\\ to
Estimate Higher Order Derivatives}

\author{David M. Bradley}
\address{Department of Mathematics and Statistics\\
         University of Maine\\
         5752 Neville Hall\\
         Orono, Maine 04469--5752\\
         U.S.A.}
\email{bradley@math.umaine.edu, dbradley@member.ams.org}

\date{Submitted
June 1, 1999. Revised August 17, 1999.}

\subjclass{Primary: 44A20; Secondary: 65D30, 65R10}

\keywords{Integral Transforms, Simpson's Rule, Estimating
Derivatives, Differentiating under the Integral Sign.}

\maketitle

\section{Introduction}
When doing error analysis for numerical quadrature, achieving good
uniform bounds on higher order derivatives of the integrand is
paramount. As undergraduates become increasingly adept with
programmable calculators, numerical integration schemes such as
Simpson's Rule and the Trapezoidal Rule take on a new relevance.
Although it may be a rare calculus class that dwells overmuch on
error bounds for such schemes, this may be due as much to the
perceived paucity of interesting examples for which decent error
bounds are readily achievable as to the general weakness in
algebraic skills necessary for the requisite understanding of
inequalities.  The purpose of this article, therefore, is to offer
some interesting, non-trivial examples for which the error
analysis, if not elegant, is at least simple enough to carry out
in the classroom.

\section{The Sine, Cosine and Exponential Integrals}
\label{sincsec}

The sine integral (\cite[pp.~231--232]{AS}
or~\cite[pp.~503--504]{CRC}) is given by
\begin{equation}
   \mathrm{Si}(x)=\int_0^x \frac{\sin t}{t}\,dt.
\label{sinc}
\end{equation}
The integrand of~(\ref{sinc}) is a scaled version of the sinc
function, as it is referred to by those who work on problems
involved with signal processing and reconstruction.  Thus,
$\mathrm{sinc}(t)=\sin(\pi t)/(\pi t)$ for $t\ne 0$, and it is
convenient to define $\sin 0/0=1$, thereby removing the removable
discontinuity.  For Simpson's rule, we require a good uniform
upper bound on the absolute value of the fourth derivative of the
integrand (\cite[p.~886]{AS} or~\cite[pp.~57--58]{DR}). Repeated
application of the product rule for differentiation gives
\[
   \D{t}{4} \frac{\sin t}{t} = \frac{\sin t}{t} + \frac{4\cos
   t}{t^2} - \frac{12\sin t}{t^3} - \frac{24\cos
   t}{t^4}+\frac{24\sin t}{t^5},
\]
which in this form is difficult to estimate, not only because of
the complicated nature of the formula, but also because of the
apparent difficulties when $t$ is near zero. However,
\begin{equation}
   \frac{\sin t}{t} = \iu \cos(st)\,ds,
\label{sincrepcos}
\end{equation}
and hence by Leibniz's rule for differentiating under the integral
sign,
\begin{equation}
   \D{t}{4}\frac{\sin t}{t} = \iu s^4\cos(st)\,ds.
\label{d4sinc}
\end{equation}
Since $|\cos(st)|\le 1$ for all real $s$ and $t$, it follows that
\begin{equation}
   \left|\D{t}{4}\frac{\sin t}{t}\right| \le \iu s^4 \,ds = \frac15.
\label{d4sincbound}
\end{equation}
The inequality is sharp, as can be seen by substituting $t=0$
in~(\ref{d4sinc}).  In a similar manner, one can show that for all
nonnegative integers $k$,
\begin{equation}
   \Bigg|\D{t}{k}\frac{\sin t}{t}\Bigg|
   \le \iu s^k \,ds = \frac{1}{k+1},
\label{dksincbound}
\end{equation}
with equality at $t=0$ when $k$ is even.  Instead of starting
with~(\ref{sincrepcos}), an analysis based on the equivalent
representation
\begin{equation}
   \frac{\sin t}{t} = \frac12\int_{-1}^1 e^{ist}\,ds
\label{sincrepexp}
\end{equation}
can likewise be given, although in freshman/sophomore calculus
classes it is less likely that students will be comfortable with
the complex exponential.

For the cosine integrals (\cite[pp.~231--232]{AS}
or~\cite[pp.~503--504]{CRC})
\begin{equation}
   \Cin(x) = \int_0^x \frac{1-\cos t}{t}\,dt,\qquad
   \mathrm{Ci}(x) = \gamma +\log x + \int_0^x \frac{\cos
   t-1}{t}\,dt,
\label{Cidef}
\end{equation}
we define the integrand to be zero when $t=0$ so that for all real
$t$,
\begin{equation}
   \frac{1-\cos t}{t}=\iu \sin(st)\,ds,
\label{cinrepsin}
\end{equation}
and hence
\[
   \D{t}{4}\frac{1-\cos t}{t} = \iu s^4 \sin(st)\,ds.
\]
Since $|\sin(st)|\le 1$ for all real $s$ and $t$, it follows that
\begin{equation}
   \left|\D{t}{4}\frac{1-\cos t}{t}\right|
   \le \iu s^4 \,ds = \frac{1}{5},
\label{d4cinbound}
\end{equation}
and in general, for all nonnegative integers $k$,
\begin{equation}
   \left|\D{t}{k}\frac{1-\cos t}{t}\right|
   \le \iu s^k \,ds = \frac{1}{k+1},
\label{dkcinbound}
\end{equation}
with equality at $t=0$ when $k$ is odd.

The hyperbolic sine and cosine integrals (\cite[p.~231]{AS}
or~\cite[p.~936]{GR})
\begin{eqnarray*}
   \mathrm{Shi}(x) &=& \int_0^x \frac{\sinh t}{t}\,dt,\\
   \mathrm{Cinh}(x) &=& \int_0^x \frac{1-\cosh t}{t}\,dt,\\
   \mathrm{Chi}(x) &=& \gamma+\log x
      +\int_0^x\frac{\cosh t-1}{t}\,dt,\qquad x>0
\end{eqnarray*}
can be treated in much the same way.  For example, for
$\mathrm{Shi}$ one replaces $\sin$ by $\sinh$ and $\cos$ by
$\cosh$ in~(\ref{sinc})--(\ref{d4sinc}). The bound
in~(\ref{dksincbound}) is modified by the presence of an extra
factor of $\cosh t$ or $|\sinh t|$ according to whether $k$ is
even or odd.  In either case, equality obtains when $t=0$.

The exponential integral (\cite[pp.~228--231]{AS}
or~\cite[pp.~504--505]{CRC})
\begin{equation}
   \mathrm{E}_1(x) = \int_1^{\infty} \frac{e^{-xt}}{t}\,dt,
   \qquad x>0
\label{E1}
\end{equation}
is another natural choice to study using this technique. In this
case, it is easier to deal with the complementary expression
\begin{equation}
   \mathrm{Ein}(x) = \int_0^x \frac{1-e^{-t}}{t}\,dt,
\label{Ein}
\end{equation}
since the latter defines an entire function. In view of the
relationship (\cite[p.~228]{AS} or~\cite[p.~40]{Olv})
\[
   \mathrm{Ein}(x)=\log x+\gamma+\mathrm{E}_1(x),
   \qquad x>0,
\]
there is no essential difference between the two. If we define the
integrand of~(\ref{Ein}) to be zero when $t=0$, then
\begin{equation}
   \frac{1-e^{-t}}{t}=\iu e^{-st}\,ds,
   \qquad \D{t}{k}\frac{1-e^{-t}}{t}=\iu (-s)^k e^{-st}\,ds,
\label{einrepexp}
\end{equation}
and so for $t\ge 0$ and $k$ a nonnegative integer,
\begin{equation}
   \left|\D{t}{k}\frac{1-e^{-t}}{t}\right|
   =\iu s^k e^{-st}\,ds
   \le \iu s^k\,ds=\frac{1}{k+1},
\label{dkeinbound}
\end{equation}
with equality again when $t=0$.

\section{How Good Are These Estimates in Practice?}
The estimates~(\ref{d4sincbound}) and~(\ref{dkeinbound}) of the
previous section are best possible in the sense that equality
holds in each when $t=0$.  On the other hand, the corresponding
derivatives~(\ref{d4sinc}) and~(\ref{einrepexp}) each tend to zero
as $t$ grows without bound, so it is clear that a uniform
numerical bound for the entire range of $t$ values is less than
ideal. Nevertheless, it is instructive to see just how well our
estimates hold up in practice.  We confine ourselves here to a
single example, the cosine integral, in which we test the
inequality~(\ref{d4cinbound}) used to give an error estimate for
Simpson's Rule against the error arising from an actual
computation.

Let $f:\R\to\R$ be given by $f(0)=0$ and $f(t)=(1-\cos t)/t$ for
$t\ne0$.  Then the cosine integral~(\ref{Cidef}) is given by
\begin{equation}
   \Cin(x)=\int_0^x f(t)\,dt.
\label{Cindef}
\end{equation}
For even positive integers $n$ and real $x>0$, define
\[
   S_n(x) :=
   \frac{x}{3n}\bigg\{f(0)+4\sum_{j=1}^{n/2}f\left(\frac{(2j-1)x}{n}\right)
   +2\sum_{j=1}^{n/2-1}f\left(\frac{2jx}{n}\right)+f(x)\bigg\},
\]
the approximation to the integral~(\ref{Cindef}) obtained by
applying Simpson's Rule with $n$ subdivisions of the interval
$[0,x]$.  The error is given by~\cite[p.~58]{DR}
\[
   E_n(x) := S_n(x)-\Cin(x) = \frac{x^5f^{(4)}(\xi)}{180 n^4},
\]
where $\xi=\xi(x)\in[0,x]$.  In view of~(\ref{d4cinbound}), we
have the inequality
\begin{equation}
   |E_n(x)| \le B_n(x),\quad
   \mbox{where}\quad B_n(x) := \frac{x^5}{900 n^4}.
\label{SimpCinBound}
\end{equation}
To compare the estimate $B_n$ with the actual error $E_n$, we let
\[
   R_n(x) := \frac{B_n(x)}{E_n(x)} =
   \frac{B_n(x)}{S_n(x)-\Cin(x)}
\]
be the ratio of the two quantities.  Then $R_n$ is defined except
for those points at which the actual error vanishes, and elsewhere
$|R_n|\ge 1$.

Using Maple's built-in cosine integral, it is possible, in
principle, to compute $E_n$ (and hence $R_n$) to any desired
precision.  Table~\ref{Tab:R} gives an indication of the range of
values of $R_n(x)$ for $0\le x\le 10$.  Tabular values were
computed using Maple V Release 5 with a conservative working
precision of twenty decimal places.  Entries were then rounded to
two decimal places to make comparison easier.  Thus, for example,
the error bound given by~(\ref{SimpCinBound}) arising from the
estimate~(\ref{d4cinbound}) exceeds the actual error by a factor
of almost $3$ when $x=5$ and Simpson's Rule is used with $n=10$,
$100$, or $1000$ subdivisions.  When $x=7$, the corresponding
overestimate is by a factor of less than $12$.

\begin{table}[h]\begin{center}
\begin{tabular}{|c|c|c|c|}
\hline
  $x \setminus n$ & 10 & 100 & 1000 \\ \hline
  1.0 & 2.55 & 2.56 & 2.56 \\ \hline
  2.0 & 1.54 & 1.55 & 1.55 \\ \hline
  3.0 & 1.43 & 1.44 & 1.44 \\ \hline
  4.0 & 1.75 & 1.77 & 1.77 \\ \hline
  5.0 & 2.77 & 2.82 & 2.82 \\ \hline
  6.0 & 5.63 & 5.75 & 5.75 \\ \hline
  7.0 & 11.31 & 11.58 & 11.58 \\ \hline
  8.0 & 10.66 & 11.18 & 11.18 \\ \hline
  9.0 & 6.97 & 7.51 & 7.51 \\ \hline
  10.0 & 5.59 & 6.14 & 6.15 \\ \hline
\end{tabular}
\bigskip
\caption{Values of $R_n(x)$} \label{Tab:R}\end{center}
\end{table}

\begin{figure}
   \begin{flushleft}
   \rotatebox{270}{\scalebox{0.55}{\includegraphics{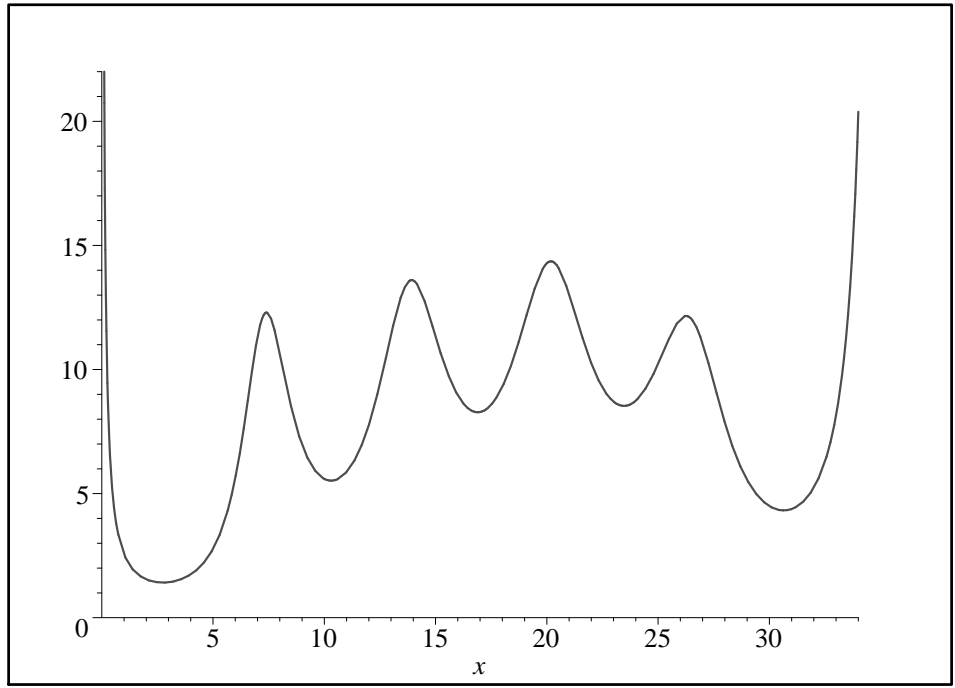}}}
   \caption{Graph of $R_{10}(x)$}
   \label{Fig:PlotR10}
   \end{flushleft}
\end{figure}

Figure~\ref{Fig:PlotR10} is a graph of $R_{10}(x)$ for $0\le x\le
34$, and makes it clear that there is nothing special about our
choice of integer values for $x$ in Table~\ref{Tab:R}. We see that
for $0\le x\le 20$, $1\le R_{10}(x)\le 15$, i.e.\ the
bound~(\ref{SimpCinBound}) is never farther than a factor of $15$
from the truth on the interval $[0,20]$.  Beyond $x=30$, we see
what appears to be the formation of a spike in the graph. This is
best explained by the computed values
$E_{10}(34.858)=0.00016504...$ and
$E_{10}(34.859)=-0.000095463...$.  Since the actual error $E_{10}$
is clearly a continuous function, the intermediate value theorem
implies that $E_{10}$ vanishes at some point in the interval
$[34.858,34.859]$. On the other hand, the estimated error $B_{10}$
is clearly an increasing function, so it follows that there exist
exceptional tiny intervals on which the estimated error
overestimates the actual error by arbitrarily large factors.
However, this kind of phenomenon is in general to be expected for
Simpson's Rule no matter what sort of error estimate is used
(apart from one based on the actual error function itself!)

\section{Understanding the Examples}
\label{tay} Although the reasoning of Section~\ref{sincsec} may be
easy to follow, the student might legitimately wonder how one
arrives at useful integral representations such as
(\ref{sincrepcos}), (\ref{sincrepexp}), (\ref{cinrepsin}), or
(\ref{einrepexp}) in the first place.  We begin by observing that
the examples of Section~\ref{sincsec} deal with integrands of the
form $(f(t)-f(0))/t$.  More generally, if $f:[a,b]\to\R$ has $n+1$
continuous derivatives, repeated integration by parts applied to
the equation
\[
   f(t)=f(a)+\int_a^t f'(u)\,du
\]
yields Taylor's formula with the integral form of the remainder:
\begin{equation}
   f(t)=\sum_{j=0}^n \frac{(t-a)^j}{j!}f^{(j)}(a)
    + \int_a^t \frac{(t-u)^n}{n!} f^{(n+1)}(u) \,du,\qquad
   a<t<b.
\label{Tay}
\end{equation}
If we make the change of variable $u=(t-a)s+a$ and
rearrange~(\ref{Tay}), we obtain
\begin{equation}
   \frac{f(t)-\sum_{j=0}^n (t-a)^j f^{(j)}(a)/j!}{(t-a)^{n+1}}
   =\iu \frac{(1-s)^n}{n!} f^{(n+1)}((t-a)s+a)\,ds,
\label{Tayform}
\end{equation}
and hence if $f$ is sufficiently differentiable, Leibniz's rule
for differentiating under the integral sign yields
\begin{eqnarray}
   \lefteqn{\D{t}{k}\frac{f(t)-\sum_{j=0}^n (t-a)^j
   f^{(j)}(a)/j!}{(t-a)^{n+1}}} \hspace{1in} \nonumber\\
&&   =\frac{1}{n!}\iu s^k (1-s)^n f^{(n+k+1)}((t-a)s+a)\,ds.
\label{general}
\end{eqnarray}
If a uniform upper bound on the absolute value of the derivative
of $f$ of order $n+k+1$ is known, say
\begin{equation}
   \sup_{a\le t\le b}\Big| f^{(n+k+1)} (t)\Big| \le M,
\label{M}
\end{equation}
then we have the upper bound
\begin{equation}
   \frac{M}{n!}\iu s^k (1-s)^n\,ds = \frac{k!\,M}{(n+k+1)!}
\label{genbound}
\end{equation}
for the absolute value of the left hand side of~(\ref{general}),
i.e.,
\begin{equation}
   \left| \D{t}{k}\frac{f(t)-\sum_{j=0}^n (t-a)^j
   f^{(j)}(a)/j!}{(t-a)^{n+1}}\right|
   \le \frac{k!\,M}{(n+k+1)!}.
\label{generalbound}
\end{equation}
It is readily apparent
that~(\ref{dksincbound}),~(\ref{dkcinbound}),
and~(\ref{dkeinbound}) are special cases of~(\ref{generalbound})
in which $a=n=0$ and $f(t)=\sin t$, $f(t)= -\cos t,$ and
$f(t)=-\exp(-t)$, respectively.

For an example with $n>0$, consider
\[
   \int_0^x \frac{1-\cos t}{t^2}\,dt,
\]
which arises in the study of incomplete versions of the Frullani
integral (\cite[p.~470]{GR} or~\cite[p.~398]{CRC})
\begin{eqnarray*}
   \il \frac{\cos \alpha t-\cos \beta t}{t^2}\,dt
   &=&\il\left( \frac{1-\cos \beta t}{t^2}-\frac{1-\cos \alpha
   t}{t^2}\right)\,dt\\
   &=& (|\beta|-|\alpha|)\frac{\pi}2,
   \qquad \alpha,\beta\in\R.
\end{eqnarray*}
Repeated application of the product rule for differentiation gives
\[
   \D{t}{4}\frac{1-\cos t}{t^2}
   =\frac{120(1-\cos t)}{t^6} -\frac{96\sin t}{t^5}
   +\frac{36\cos t}{t^4}+ \frac{8\sin t}{t^3}
   -\frac{\cos t}{t^2},
\]
which again, is troublesome to estimate as it stands. On the other
hand,~(\ref{general}) with $f(t)=\cos t$, $n=1$, $a=0$, and $k=4$
gives
\[
   \D{t}{4}\frac{1-\cos t}{t^2} = \iu s^4 (1-s)\cos(st)\,ds,
\]
whence
\[
   \left|\D{t}{4}\frac{1-\cos t}{t^2}\right|
   \le \iu s^4 (1-s)\,ds = \frac1{30},
\]
with equality at $t=0$.  More generally,~(\ref{M}) with $M=1$ is
satisfied when $f(t)=\cos t$, and so~(\ref{genbound})
and~(\ref{generalbound}) give
\[
\left|\D{t}{k}\frac{1-\cos t}{t^2}\right|
   \le \iu s^k (1-s)\,ds = \frac1{(k+1)(k+2)},
\]
with equality at $t=0$ when $k$ is even.

\section{Some More Sophisticated Examples}
The inverse tangent integral
\[
   \mathrm{Ti}_2(x) = \int_0^x \frac{\tan^{-1}(t)}{t}\,dt,
\]
would appear to be another example amenable to our technique.
Lewin~\cite[chapter~2]{Lew} devotes a whole chapter to the study
of its properties.  The inverse tangent integral evidently also
attracted the interest of Ramanujan~\cite{Rama}, who among other
things used it to develop a rapidly convergent series of
hyperbolic functions for Catalan's constant, $\mathrm{Ti}_2(1)$
(\cite[p.~807]{AS},~\cite{DMB1}). Formula~(\ref{Tayform}) with
$a=n=0$ and $f(t)=\tan^{-1}(t)$ yields
\begin{equation}
   \frac{\tan^{-1}(t)}{t} = \iu \frac{ds}{1+s^2 t^2}.
\label{atn}
\end{equation}
After differentiating under the integral and simplifying the
resulting expression, one arrives at
\begin{equation}
   \D{t}{4}\frac{\tan^{-1}(t)}{t}
   = \iu \frac{24 s^4(5 s^4 t^4- 10 s^2 t^2+ 1)}{(1+s^2 t^2)^5}\,ds.
\label{d4ti}
\end{equation}
It is difficult to obtain a good uniform estimate on the size of
the integrand in~(\ref{d4ti}), so we return to~(\ref{atn}). After
making the change of variable $s=1/u$ and employing the formula
\[
   \il e^{-uv}\cos(vt)\,dv = \frac{u}{u^2+t^2},\qquad u>0,
\]
an application of Fubini's theorem~\cite[p.~67]{Folland} gives
\begin{eqnarray*}
   \frac{\tan^{-1}(t)}{t}
   &=& \int_1^\infty \frac{1}{u}\cdot
   \frac{u\,du}{u^2+t^2}\\
   &=&\int_1^\infty \frac1{u}\il
   e^{-uv}\cos(vt)\,dv\,du\\
   &=&\il\int_1^\infty
   \frac{e^{-uv}}{u}\,du\cos(vt)\,dv\\
   &=&\il \mathrm{E}_1(v)\cos(vt)\,dv,
\end{eqnarray*}
where $\mathrm{E}_1$ is the exponential integral~(\ref{E1}). Thus,
\[
   \D{t}{4}\frac{\tan^{-1}(t)}{t}
   = \il v^4\, \mathrm{E}_1(v)\cos(vt)\,dv,
\]
and hence by Tonnelli's theorem~\cite[p.~67]{Folland},
\begin{eqnarray*}
   \left|\D{t}{4}\frac{\tan^{-1}(t)}{t}\right|
   &\le& \il v^4\,\mathrm{E}_1(v)\,dv\\
   &=&\il v^4\int_1^\infty \frac{e^{-uv}}{u}\,du\,dv\\
   &=&\int_1^\infty\frac1{u}\il v^4 e^{-uv}\,dv\,du\\
   &=&4!\int_1^\infty\frac{du}{u^6}\\
   &=&\frac{24}5.
\end{eqnarray*}
In the same manner, it can be shown more generally that for all
nonnegative integers $k$,
\[
   \left|\D{t}{k}\frac{\tan^{-1}(t)}{t}\right|
   \le \il v^k\, \mathrm{E}_1(v)\,dv = \frac{k!}{k+1},
\]
with equality at $t=0$ when $k$ is even.

The previous example gives an idea of what is possible when one
steps outside the framework of Section~\ref{tay}.  As another
example, suppose one wanted to estimate the even derivatives of
the tangent function.  From~\cite[p.~388]{GR},
\[
   \tan t = 2\il \frac{\sinh 2st}{\sinh \pi s}\,ds,\qquad
   -\pi/2<t<\pi/2.
\]
It follows that for all nonnegative integers $k$, we have
\begin{equation}
   \D{t}{2k} \tan t=2\il (2s)^{2k} \frac{\sinh 2st}{\sinh \pi
   s}\,ds, \qquad -\pi/2<t<\pi/2.
\label{tan2k}
\end{equation}
As $s$ gets large, we expect
\[
   \frac{\sinh 2st}{\sinh \pi s}
   =\frac{e^{2st}-e^{-2st}}{e^{\pi s}-e^{-\pi s}}
\]
to behave like $e^{(2t-\pi)s}$ for $0< 2t<\pi$.  In fact, it is
easy to prove that the inequality
\[
   \frac{\sinh 2s|t|}{\sinh \pi s}
   \le e^{(2|t|-\pi)s}
\]
holds for $-\pi<2t<\pi$ and $s>0.$ Therefore, by~(\ref{tan2k}), we
have
\begin{equation}
   \left|\D{t}{2k} \tan t\right|
   \le 2\il (2s)^{2k} e^{(2|t|-\pi)s}\,ds
   = \frac{(2k)!}{(\pi/2-|t|)^{2k+1}}
\label{tan2kest}
\end{equation}
for all nonnegative integers $k$. The estimate~(\ref{tan2kest}) is
remarkably good even when $|t|$ is small, with asymptotic equality
in the limit as $|t|\to\pi/2-$.

The previous two examples succeeded because we were able to
represent the desired function as a Fourier (respectively,
Laplace) transform of a well-behaved function.   As a final
example, consider the integral
\begin{equation}
   \int_a^b t^{-2\kappa}e^{-ut}e^{\kappa\, {\mathrm{Ein}}(t)}\,dt,
\label{qint}
\end{equation}
in which $0<a<b$, $\kappa>1$, and $u>0$ are real parameters, and
Ein is the complementary exponential integral~(\ref{Ein}).  The
integral~(\ref{qint}) arises in the solution of certain advanced
argument difference-differential equations relating to
sieves~(\cite{DMB2},~\cite{DMB3}), and as such, it is desirable to
be able to compute it accurately for various values of the
parameters.  To study the integrand of~(\ref{qint}), we let $c>0$
and define
\begin{equation}
   \lambda_{\kappa}(v) := \frac{1}{2\pi i}
   \int_{c-i\infty}^{c+i\infty} e^{vz}
   z^{-2\kappa}e^{\kappa\, {\mathrm{Ein}}(z)}\,dz.
\label{lineint}
\end{equation}
By Cauchy's theorem, the line integral~(\ref{lineint}) is
independent of $c$, and vanishes for $v\le 0$.  It is not hard to
show that $\lambda_{\kappa}$ is continuous and satisfies the delay
differential equation
\begin{equation}
   \left(v^{1-\kappa}\lambda_{\kappa}(v)\right)' =\kappa v^{-\kappa}
   \lambda_{\kappa}(v-1),\qquad v>0,
\label{delay}
\end{equation}
with boundary condition
\begin{equation}
   \lambda_{\kappa}(v)
   = \frac{e^{\kappa\gamma}}{\Gamma(\kappa)} v^{\kappa-1},
   \qquad 0\le v\le 1.
\label{boundary}
\end{equation}
Let $f(t)$ denote the integrand of~(\ref{qint}). By Laplace
inversion,
\[
   \il e^{-vt}\lambda_{\kappa}(v)\,dv =
   t^{-2\kappa}e^{\kappa\,{\mathrm{Ein}}(t)},\qquad t>0,
\]
and differentiating under the integral sign to obtain
\begin{eqnarray*}
   f^{(4)}(t)&=&\D{t}4 t^{-2\kappa}e^{-ut}e^{\kappa\,{\mathrm{Ein}}(t)}
   = \il (u+v)^4 e^{-(u+v)t}\lambda_{\kappa}(v)\,dv,\\
   f^{(5)}(t)&=&\D{t}5t^{-2\kappa}e^{-ut}e^{\kappa\,{\mathrm{Ein}}(t)}
   =- \il (u+v)^5 e^{-(u+v)t}\lambda_{\kappa}(v)\,dv,
\end{eqnarray*}
can be justified.  But, the delay differential
equation~(\ref{delay}) and boundary condition~(\ref{boundary})
together show that $\lambda_{\kappa}$ is nonnegative, and hence
$f^{(4)}$ is nonnegative and nonincreasing. It follows that
\[
   \sup_{a\le t\le b} \left|f^{(4)}(t)\right| = f^{(4)}(a),
\]
a considerable simplification.

\section{Conclusion}
Of course, there are other methods for bounding the derivatives of
a suitable function, the most familiar of which is undoubtedly
Cauchy's inequality (see eg.\ \cite[p.~91]{WW}) from the theory of
complex variables.  One can also consider alternative approaches
to error analysis which do not involve estimating higher order
derivatives. Chapter 4 of Davis and Rabinowitz~\cite{DR} is
devoted to error analysis for various approximate integration
schemes.  In addition to a section on error estimates via analytic
function theory, alluded to previously, there is a lovely section
(see pp.~317--332) describing applications of functional analysis
to numerical integration and error estimation.

In this paper, however, we have deliberately focused on the use of
integral transforms arising primarily in the context of real
variable theory.  The technique can be summarized as follows: find
a suitable integral representation for the function whose
derivatives are to be estimated, differentiate under the integral
sign, and estimate the resulting integral.  The technique is
hardly new, but is seldom used to its fullest advantage.  It is
hoped that some of the examples provided herein could be used to
enrich the discussion of numerical integration in a typical
calculus or numerical analysis course.

\end{document}